\documentclass[a4paper,oneside,reqno]{amsart}
\usepackage{xy}
\xyoption{all}
\title[Structure Constants of the Schur Algebra]{Structure Constants\\of the Schur Algebra}
\newcommand{\ii}{\mathbf i}

\DeclareMathOperator{\End}{End}
\newtheorem{proposition}{Proposition}
\theoremstyle{remark}
\newtheorem{example}{Example}
\usepackage[margin=3.5cm]{geometry}
\title{Graphic Interpretation of the Structure Constants of the Schur Algebra}
\author{T. Geetha}
\author{Amritanshu Prasad}
\begin{document}
 ICWM Poster Abstract Number: 20140020
\maketitle
\begin{center}
 The Institute of Mathematical Sciences\\
 Chennai, India\\
 E-mail: geetha\textunderscore curie@yahoo.co.in\\
  amri@imsc.res.in
\end{center}

\section*{The Schur Algebra}

Let $k$ be a field, and let $V$ denote the vector space $k^n$.
Let $e_1,\dotsc, e_n$ denote the coordinate vectors of $V$.
Let $V^{\otimes d}$ denote the $d$-fold tensor product of $V$ with itself.
Let
\begin{equation*}
  I(n,d) = \{\ii = (i_1,\dotsc,i_d)\mid 1\leq i_s\leq n \text{ for } 1\leq s \leq d\}.
\end{equation*}
Let $e_\ii = e_{i_1}\otimes \dotsc \otimes e_{i_d}$ for each $\ii = (i_1,\dotsc, i_d)\in I(n,d)$.
Then the vectors $e_\ii$, as $\ii$ runs over $I(n,d)$ form a basis for $V^{\otimes d}$.
The symmetric group $S_d$ acts on $I(n,d)$ by
\begin{equation*}
  w\cdot (i_1,\dotsc, i_d) = (i_{w(1)}, \dotsc, i_{w(d)}).
\end{equation*}
Since $S_d$ permutes the elements of a basis of $V^{\otimes d}$, it becomes a permutation representation $k[I(n,d)]$ of $S_d$.

For our purposes, define the Schur algebra $S(n,d)$ by:
\begin{equation*}
  S(n,d) = \End_{S_d} V^{\otimes d} = \End_{S_d} k[I(n,d)].
\end{equation*}
For each $S_d$-orbit $O\in I(n,d)^2$, and $f\in k[I(n,d)]$, define
\begin{equation*}
  T_O f(x) = \sum_{\{y\mid (x,y)\in O\}} f(y).
\end{equation*}
Then $T_O$, as $O$ runs over the set of $S_d$-orbits in $I(n,d)^2$, is a basis for $S(n,d)$.
For any $S_d$-orbit $O$ in $I(n,d)^2$, take any $(x,y)\in O$.
Then the coefficient of $T_O$ in the expansion of $T_{O_1}\circ T_{O_2}(x, y)$ is
\begin{equation*}
  \#\{z \in I(n,d)\mid (x, z)\in O_1,\; (z, y)\in O_2\}.
\end{equation*}

\section*{The balls-in-boxes picture}

Let $B(n,d)$ denote the set of all possible ways of putting $d$ distinguishable balls (say, numbered $1, \dotsc, d$) in $n$ boxes numbered $1,\dotsc, n$.

Given $\ii = (i_1,\dotsc, i_d)\in I(n,d)$ define $\phi(\ii) \in B(n,d)$ by putting the $s$th ball in the $i_s$th box for each $1\leq s\leq d$.
The index $i_s$ can be recovered from $\phi(\ii)$ as the number of the box in which the $s$th ball is.
A balls-in-boxes configuration of $d$ balls in $n$ boxes is denoted by a word with $n+1$ occurences of the symbol \lq$|$\rq{} and with one occurence of each of the integers $1,\dotsc,d$ such that
\begin{enumerate}
\item the word begins and ends with \lq$|$\rq{}.
\item the integers between two \lq$|$\rq{}s are in increasing order.
\end{enumerate}
The $n+1$ \lq$|$\rq{}s form the boundaries of the $n$ boxes.
The integers that appear between the $i$th and $(i+1)$st \lq$|$\rq{} are the balls which are in the $i$th box.

\begin{example}
  Consider $n=2$ and $d=4$.
  Then
  \begin{equation*}
    \begin{array}{|c|c|}
      \hline
      \ii & \text{balls-in-boxes}\\
      \hline
      (1, 1, 1, 1) & |1234||\\
      (1, 1, 1, 2) & |123|4|\\
      
      \hline
    \end{array}
  \end{equation*}
\end{example}

The symmetric group $S_d$ acts on $B(n,d)$ by simply permuting the $d$ balls.
\begin{proposition}
  The function $\phi:I(n,d)\to B(n,d)$ is an isomorphism of $S_d$-sets.
\end{proposition}
We may therefore think of the Schur algebra $S(n,d)$ as the endomorphism algebra of the permutation representation of $S_d$ coming from its action on $B(n,d)$:
\begin{equation*}
  S(n,d) = \End_{S_d} k[B(n,d)].
\end{equation*}

As a consequence, $S(n,d)$ has a basis indexed by orbits for the diagonal action of $S_d$ on $B(n,d)^2$.
Given $(a,b)\in B(n,d)^2$, define $\Gamma(a,b)$ to be the bipartite graph with $2n$ vertices arranged in two rows of $n$ vertices each (the vertices in the upper row correspond to the boxes for $b$ and the vertices in the lower row correspond to the boxes for $a$) and $d$ edges (one for each ball).
If the $s$th ball lies in the $i$th box for $a$ and the $j$th box for $b$, then the edge corresponding to the $s$th ball joins the $j$th vertex in the first row to the $i$th vertex in the second row.
Note that the valency of the $i$th edge in the bottom row of $\Gamma_{a,b}$ is the number of balls in the $i$th box for $a$, and the valency of the $j$th edge in the top row of $\Gamma_{a,b}$ is the number of balls in the $j$th box for $b$ for all $i,j\in \{1, \dotsc, n\}$.
\begin{example}
  If $a = |123|4|$ and $b = |13|24|$, then the graph $\Gamma_{a,b}$ is given by
 \small
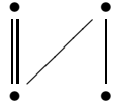
{ \begin{equation*}
    \xymatrix{
      \bullet \ar@{=}[d] & \bullet\ar@{-}[dl] \ar@{-}[d]\\
      \bullet & \bullet
    }
  \end{equation*}}
\end{example}
Since the action of $S_d$ permutes the $d$ balls, $(a,b)$ and $(a',b')$ in $B(n,d)^2$ lie in the same $S_d$-orbit if $(a',b')$ can be obtained from $(a,b)$ by renaming the balls.
The following proposition easily follows:
\begin{proposition}
  Let $\Gamma(n,d)$ denote the set of all bipartite graphs with $n$ vertices in each part and $d$ edges.
  The map $\gamma:(a,b)\mapsto \Gamma_{a,b}$ descends to a bijection
  \begin{equation*}
    \bar\gamma:S_d\backslash B(n,d)^2 \to \Gamma(n,d).
  \end{equation*}
\end{proposition}
For $\Gamma \in \Gamma(n,d)$ define $\xi_\Gamma\in S(n,d)$ by
\begin{equation*}
  \xi_\Gamma f (a) = \sum_{\{b\mid \gamma(a,b)=\Gamma\}}f(b), \text{ for each } f\in k[B(n,d)].
\end{equation*}
By the general theory of permutation representations, the set
\begin{equation*}
  \{\xi_\Gamma \mid \Gamma \in \Gamma(n,d)\}
\end{equation*}
forms a basis for $S(n,d)$.

For each $a\in B(n,d)$, let $e_a\in k[B(n,d)]$ denote the function which is $1$ at $a$ and $0$ everywhere else.
Then
\begin{eqnarray*}
  \xi_\Gamma e_b(x) & = \sum_{\{y\mid \gamma(x,y)=\Gamma\}} e_b(y)\\
  & =
  \begin{cases}
    1 & \text{ if } \gamma(x,b) = \Gamma\\
    0 & \text{ otherwise.}
  \end{cases}
\end{eqnarray*}
Therefore
\begin{equation}
  \label{eq:action}
  \xi_\Gamma e_b = \sum_{\{a\mid \gamma(a, b) = \Gamma\}} e_a.
\end{equation}
In particular, if the valency of $j$th vertex in the first row of $\Gamma$ does not match the number of balls in the $j$th box for $b$ for any $j\in \{1, \dotsc, n\}$, then $\xi_\Gamma e_b = 0$.

Consider another example:
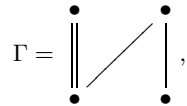
\small{\begin{equation*}
  \Gamma =     \vcenter{\xymatrix{
      \bullet \ar@{=}[d] & \bullet\ar@{-}[dl] \ar@{-}[d]\\
      \bullet & \bullet
    }},
\end{equation*}}
then
\begin{equation*}
  \xi_\Gamma e_{|12|34|} = e_{|123|4|} + e_{|124|3|}.
\end{equation*}
This rule has a nice visualization.
Replace the nodes in the first row of $\Gamma$ by the contents of the corresponding boxes for the balls-in-boxes configuration $b$:
\small
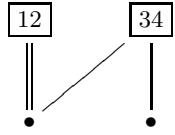
{\begin{equation*}
  \xymatrix{
      \framebox{12} \ar@{=}[d] & \framebox{34}\ar@{-}[dl] \ar@{-}[d]\\
      \bullet & \bullet
    }
\end{equation*}}
The element $\xi_\Gamma e_{|12|34|}$ is the sum of $e_a$ over all distinct configurations $a\in B(n,d)$ which can be obtained by sending each of the $d$ balls along an edge to the second row of nodes in $\Gamma$ in such a way that each of the $d$ edges is used exactly once.
In the running example, the balls numbered $1$ and $2$ have to be sent to the first node of the second row, and one of the balls numbered $3$ and $4$ can be sent to the first node, while the other is sent to the second node.
Thus there are two possible configurations:
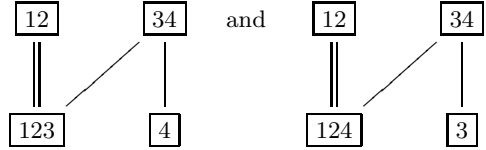
\small{\begin{equation*}
  \xymatrix{
      \framebox{12} \ar@{=}[d] & \framebox{34}\ar@{-}[dl] \ar@{-}[d]\\
      \framebox{123} & \framebox{4}
    }
    \quad\text{ and }\quad
  \xymatrix{
      \framebox{12} \ar@{=}[d] & \framebox{34}\ar@{-}[dl] \ar@{-}[d]\\
      \framebox{124} & \framebox{3}
    }
\end{equation*}}
whence
\begin{equation*}
  \xi_\Gamma e_{|12|34|} = e_{|123|4|} + e_{|124|3|}.
\end{equation*}

We are now ready to describe the structure constants for $S(n,d)$.
Given two graphs $\Gamma_1$ and $\Gamma_2$ in $\Gamma(n,d)$, suppose that
\begin{equation*}
  \xi_{\Gamma_1}\xi_{\Gamma_2} = \sum_{\Gamma}c^{\Gamma}_{\Gamma_1 \Gamma_2}\xi_\Gamma.
\end{equation*}
Take any $c\in B(n,d)$ such that the $j$th box for $c$ has as many balls as the valency of the first row of $\Gamma_2$.
Take $a\in B(n,d)$ to be any configuration such that $\gamma(a,c) = \Gamma$.
Then $\xi_\Gamma e_c(a) = 1$.
On the other hand, by (\ref{eq:action}), $\xi_{\Gamma_1}\xi_{\Gamma_2} e_c$ is equal to the cardinality of the set
\begin{equation*}
  \{ b\in B(n,d) \mid \gamma(a,b) = \Gamma_1, \; \gamma(b,c) = \Gamma_2\},
\end{equation*}
whence we have
\begin{proposition}
  Given graphs $\Gamma_1, \Gamma_2, \Gamma\in \Gamma(n,d)$, the coefficient $c^\Gamma_{\Gamma_1\Gamma_2}$ of $\xi_\Gamma$ in $\xi_{\Gamma_1}\xi_{\Gamma_2}$ is given by
  \begin{equation*}
    c^\Gamma_{\Gamma_1\Gamma_2} = \# \{ b\in B(n,d) \mid \gamma(a,b) = \Gamma_1, \; \gamma(b,c) = \Gamma_2\},
  \end{equation*}
  where $(a,c)\in B(n,d)^2$ is any pair such that $\gamma(a,c) = \Gamma$.
\end{proposition}
This rule for computing the structure constants of the Schur algebra has a nice visualization.
For example, consider
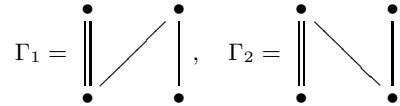
\small{\begin{equation*}
  \Gamma_1 = \vcenter{
    \xymatrix{
      \bullet \ar@{=}[d] & \bullet\ar@{-}[dl] \ar@{-}[d]\\
      \bullet & \bullet
    }
  },\quad
  \Gamma_2 = \vcenter{
    \xymatrix{
      \bullet \ar@{=}[d] \ar@{-}[dr] & \bullet \ar@{-}[d]\\
      \bullet & \bullet
    }
  }
\end{equation*}}
We may always take $c = |123|4|$.
To find the coefficient of $\xi_\Gamma$ with
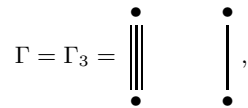
\small{\begin{equation*}
  \Gamma = \Gamma_3 = \vcenter{
    \xymatrix{
      \bullet \ar@3{-}[d] & \bullet \ar@{-}[d]\\
      \bullet & \bullet
    }
  },
\end{equation*}}
we may take $a = |123|4|$ as well.
This coefficient is the number of ways of filling in the boxes in the middle row of the diagram
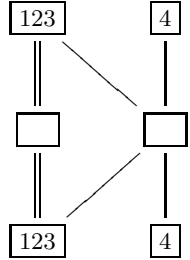
\small{\begin{equation*}
  \xymatrix{
    \framebox{123} \ar@{=}[d] \ar@{-}[dr] & \framebox{4} \ar@{-}[d]\\
    \framebox{\phantom{12}} \ar@{=}[d] & \framebox{\phantom{34}}\ar@{-}[dl] \ar@{-}[d]\\
    \framebox{123} & \framebox{4}
  }
\end{equation*}}
which are compatible with the top and bottom rows.
For this there are clearly three possibilities, namely, we can choose which of the first three balls ends up in the second box of the middle row:
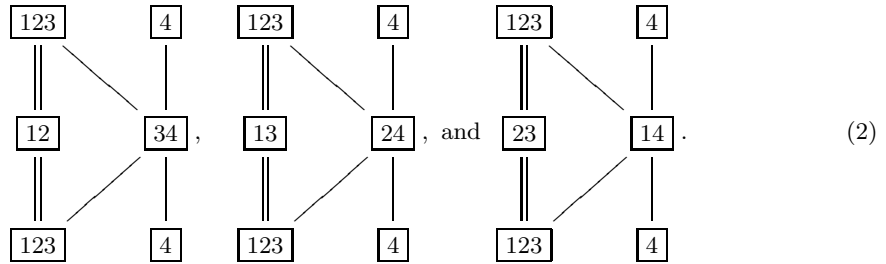
\small{\begin{equation}
  \label{eq:middle-box-3}
  \vcenter{\xymatrix{
    \framebox{123} \ar@{=}[d] \ar@{-}[dr] & \framebox{4} \ar@{-}[d]\\
    \framebox{12} \ar@{=}[d] & \framebox{34}\ar@{-}[dl] \ar@{-}[d]\\
    \framebox{123} & \framebox{4}
  }},\quad
  \vcenter{\xymatrix{
    \framebox{123} \ar@{=}[d] \ar@{-}[dr] & \framebox{4} \ar@{-}[d]\\
    \framebox{13} \ar@{=}[d] & \framebox{24}\ar@{-}[dl] \ar@{-}[d]\\
    \framebox{123} & \framebox{4}
  }}, \text{ and }
  \vcenter{\xymatrix{
    \framebox{123} \ar@{=}[d] \ar@{-}[dr] & \framebox{4} \ar@{-}[d]\\
    \framebox{23} \ar@{=}[d] & \framebox{14}\ar@{-}[dl] \ar@{-}[d]\\
    \framebox{123} & \framebox{4}
  }}.
\end{equation}}
On the other hand, for
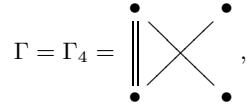
\small{\begin{equation*}
 \Gamma = \Gamma_4 =
 \vcenter{
   \xymatrix{
     \bullet \ar@{=}[d] \ar@{-}[dr] & \bullet \ar@{-}[dl]\\
     \bullet & \bullet
   }
 },
\end{equation*}}
and $c = |123|4|$, we may take $a = |124|3|$.
We need to fill in the boxes in the middle row of the diagram
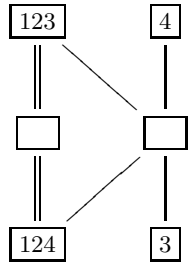
\small{\begin{equation*}
  \xymatrix{
    \framebox{123} \ar@{=}[d] \ar@{-}[dr] & \framebox{4} \ar@{-}[d]\\
    \framebox{\phantom{12}} \ar@{=}[d] & \framebox{\phantom{34}}\ar@{-}[dl] \ar@{-}[d]\\
    \framebox{124} & \framebox{3}
  }
\end{equation*}}
for which there is only one possibility, namely,
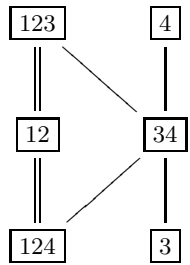
\small{\begin{equation*}
  \xymatrix{
    \framebox{123} \ar@{=}[d] \ar@{-}[dr] & \framebox{4} \ar@{-}[d]\\
    \framebox{12} \ar@{=}[d] & \framebox{34}\ar@{-}[dl] \ar@{-}[d]\\
    \framebox{124} & \framebox{3}
  }
\end{equation*}}
It turns out that for no other $\Gamma \in \Gamma(n,d)$ is it possible to find even one compatible way of filling in the middle boxes.
Indeed, if we are able to compatibly fill in the second and third rows of the diagram
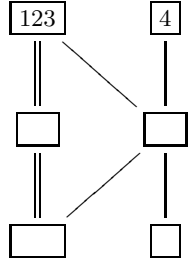
\small{\begin{equation*}
  \xymatrix{
    \framebox{123} \ar@{=}[d] \ar@{-}[dr] & \framebox{4} \ar@{-}[d]\\
    \framebox{\phantom{12}} \ar@{=}[d] & \framebox{\phantom{34}}\ar@{-}[dl] \ar@{-}[d]\\
    \framebox{\phantom{124}} & \framebox{\phantom{3}}
  }
\end{equation*}}
each ball in the first row ends up taking path consisting of an edge from $\Gamma_2$ and then an edge from $\Gamma_1$ to end up in the third row.
The path of this ball determines a bijection $\phi$ from the set of edges of $\Gamma_2$ to the set of edges of $\Gamma_1$ such that the lower node of $e$ coincides with the upper node of $\phi(e)$ for every edge $e$ of $\Gamma_2$.
Each such function (which we shall call an \emph{Euler function}) determines a graph $\Gamma$ as follows:
for each edge $e$ of $\Gamma_2$, $\Gamma$ has an edge going from the upper node of $e$ to the lower node of $\phi(e)$ (intuitively, this edge is the composition of $e$ and $\phi(e)$).
Thus all the graphs in the support of $\xi_{\Gamma_1}\xi_{\Gamma_2}$ come from an Euler function in this way.

It follows that
\begin{equation*}
  \xi_{\Gamma_1}\xi_{\Gamma_2} = 3\xi_{\Gamma_3} + \xi_{\Gamma_4}.
\end{equation*}

What we have just seen is a way to write down the structure constants of the Schur algebra.
This method, however, depends on a choice of $c\in B(n,d)$ which is compatible with the first row of $\Gamma_2$.
Of course, there is a canonical choice for $c$ (namely numbering the balls in increasing order, starting from the first box and going up to the $n$th box, as we have done in all the examples here.

In \cite{Mendez}, M\'endez gave a description for the structure constants of the Schur algebra which does not use any choice (canonical or otherwise) of $c$, but rather a labelling of $\Gamma_1$ and $\Gamma_2$, which we will now describe:

Label the edges of $\Gamma_1$ and $\Gamma_2$ in such a way that the multiple edges between the same vertices have the same label.
It does not matter is the same label is used for an arrow of $\Gamma_1$ and an arrow of $\Gamma_2$.
For instance we may use:
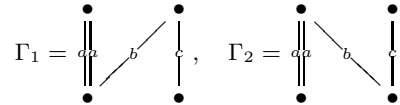
\small{\begin{equation*}
  \Gamma_1 = \vcenter{
    \xymatrix{
      \bullet \ar@{=}[d]|{aa} & \bullet\ar@{-}[dl]|b \ar@{-}[d]|c\\
      \bullet & \bullet
    }
  },\quad
  \Gamma_2 = \vcenter{
    \xymatrix{
      \bullet \ar@{=}[d]|{aa} \ar@{-}[dr]|b & \bullet \ar@{-}[d]|c\\
      \bullet & \bullet
    }
  }
\end{equation*}}
Each of the configurations in (\ref{eq:middle-box-3}) gives rise to an $n\times n$ matrix of words as follows: the word in the $(s,t)$th entry of the matrix is the sequence of ordered pairs $(e_2(i), e_1(i))$ of edge labels, $e_2(i)$ being the edge of $\Gamma_2$ and $e_1(i)$ being the edge of $\Gamma_1$ through which the $i$th ball passes to reach the top row to the bottom row, as $i$ runs over balls which begin in the $t$th box of the top row and end in the $s$th box of the bottom row taken in increasing order.
For example, corresponding to the three configurations of (\ref{eq:middle-box-3}), we get the three matrices of words:
\small{$\begin{pmatrix}
(a,a)(a,a)(b,b) & \\
 & (c,c)
\end{pmatrix}$, 
$\begin{pmatrix}
(a,a)(b,b)(a,a) & \\
 & (c,c)
\end{pmatrix}$, and
$\begin{pmatrix}
 (b,b)(a,a)(a,a) & \\
  & (c,c)
 \end{pmatrix}
$.}
On the other hand, given an $n\times n$ matrix of words, the paths taken by the balls in the $j$th box of the upper node of $\Gamma_2$ to reach the $i$th box of the lower node of $\Gamma_1$ are as follows: take the balls in the $j$th box of the upper node of $\Gamma_2$ in increasing order, then the path of each ball is given by the corresponding ordered pair in the word in the $(i,j)$th position.
This rule allows us to recover the balls-in-boxes configuration in the middle row from which this matrix of words came.
\begin{proposition}
  [M\'endez~\cite{Mendez}]
  Fix a labelling of edges of $\Gamma_1$ and $\Gamma_2$ as described above.
  The coefficient of $\xi_\Gamma$ in $\xi_{\Gamma_1}\xi_{\Gamma_2}$ is the number of $n\times n$ matrices whose $(i,j)$th entry is a word consisting of pairs $(e, f)$ which comes from an Euler function.
\end{proposition}

\end{document}